\documentclass[letterpaper, 10 pt, conference]{ieeeconf}  

\IEEEoverridecommandlockouts                              

\usepackage{cite}
\usepackage{color}
\usepackage{amsfonts,amssymb}
\usepackage{bm}
\usepackage{amsmath}
\usepackage{multirow}
\usepackage{algorithmic}
\usepackage{graphicx}
\usepackage{subfigure}
\usepackage{float}
\usepackage{textcomp}
\usepackage{graphicx} 
\usepackage{epstopdf}
\usepackage{array}
\usepackage[thmmarks,amsmath]{ntheorem}
\newtheorem{mythm}{Theorem}
\newtheorem{mydef}{Definition}
\newtheorem{assumption}{Assumption}
\newtheorem{lemma}{Lemma}

\newtheorem{proposition}{Proposition}

\usepackage{array}
\usepackage{enumerate}
\usepackage{booktabs}
\usepackage{algorithm}
\usepackage{algorithmic}
\usepackage{setspace}
\usepackage{float}
\usepackage{cases}
\usepackage{stfloats}

\usepackage{mathrsfs}

\title{\LARGE \bf
Non-convex potential games for finding\\
global solutions to sensor network localization
}

\author{Gehui Xu$^{1}$, Guanpu Chen$^{2}$, Yiguang Hong$^{3}$,  Baris Fidan$^{4}$, Thomas Parisini$^{5}$, and Karl H. Johansson$^{2}$
\thanks{*This work was supported  by 
Swedish Research Council Distinguished Professor Grant 2017-01078, Knut and Alice Wallenberg Foundation Wallenberg Scholar Grant,  Swedish Strategic Research Foundation SUCCESS Grant FUS21-0026, 
and also supported 
in part by the Digital Futures Scholar-in-Residence Program, by the European Union’s Horizon 2020 research and
innovation programme under grant agreement no. 739551 (KIOS CoE), and
by the Italian Ministry for Research in the framework of the 2017 Program
for Research Projects of National Interest (PRIN), Grant no. 2017YKXYXJ.}
\thanks{$^{1}$Gehui Xu
is with the Key Laboratory of Systems and Control, Academy of Mathematics and Systems Science, Beijing, China, and also  with the Division of Decision and Control Systems, School of	Electrical Engineering and Computer Science, KTH Royal Institute of	Technology, 100 44, Stockholm, Sweden. 
        {\tt\small xghapple@amss.ac.cn}}%
\thanks{$^{2}$ Guanpu Chen, and Karl H. Johansson 
are with the Division of Decision and Control Systems, School of	Electrical Engineering and Computer Science, KTH Royal Institute of	Technology, 100 44, Stockholm, Sweden. 
        {\tt\small  guanpu@kth.se, and kallej@kth.se}}%
        \thanks{$^{3}$Yiguang Hong is with
		Department of Control Science and Engineering, Tongji University,  Shanghai 201804, China, and is also with Shanghai Research Institute for Intelligent Autonomous Systems, Shanghai 201210 China.
        {\tt\small yghong@iss.ac.cn}}%
\thanks{$^{4}$Baris Fidan is with the Department of Mechanical and Mechatronics Engineering, University of Waterloo, Waterloo, ON N2L 3G1, Canada.
        {\tt\small fidan@uwaterloo.ca}}%
         \thanks{$^{5}$Thomas Parisini is with  the Department of Electrical and Electronic Engineering,
	Imperial College London, London SW7 2AZ, UK, and also with the Department
	of Engineering and Architecture, University of Trieste, Trieste 34127, Italy.
        {\tt\small t.parisini@imperial.ac.uk}}%
}

\begin{document}

\maketitle
\thispagestyle{empty}
\pagestyle{empty}

\begin{abstract}
Sensor network localization (SNL) problems require determining the physical coordinates of all sensors in a network. This process relies on the global coordinates of anchors and the available measurements between non-anchor and anchor nodes. Attributed to the intrinsic non-convexity, obtaining a globally optimal solution to SNL 
is challenging, as well as implementing corresponding algorithms. In this paper, we formulate a non-convex multi-player potential game for a generic SNL problem to investigate the identification condition of the global Nash equilibrium (NE) therein, where the global NE represents the global solution of SNL. We employ canonical duality theory to transform the non-convex game into a complementary dual problem. Then we develop a conjugation-based algorithm to compute the stationary points of the complementary dual problem. On this basis,
 we show an identification condition of the global NE: the stationary point of the proposed algorithm satisfies a duality relation. Finally, simulation results are provided to validate the effectiveness of the theoretical results.

\end{abstract}

\section{INTRODUCTION}

Wireless sensor networks (WSNs), due to their capabilities of sensing, processing, and communication, have a wide range of applications \cite{akyildiz2002wireless,liu2015event}, such as target tracking and detection \cite{marechal2010joint,jing2021angle}, environment monitoring \cite{sun2010corrosion}, area exploration \cite{sun2005reliable}, data collection and cooperative robot tasks \cite{jing2018weak}. {
		For all of these applications, it is essential to determine the location of every sensor with the desired accuracy. 
Estimating locations of the sensor nodes based on measurements between neighboring nodes has attracted many research interests in recent years, see typical examples \cite{estrin2000embedding,fang2023distributed}. 
}
Range-based methods constitute a common inter-node measurement approach utilizing signal transmission based techniques such as
time of arrival, time-difference of arrival, and strength of received radio frequency signals  \cite{mao2007wireless}.
Due to limited
transmission power, the measurements can only be obtained
within a radio range. A pair of nodes are called neighbors if their distance is less than this radio range \cite{wan2019sensor}. 
Also, there are some anchor nodes whose global positions are known
\cite{ahmadi2021distributed}. Then a sensor network localization
(SNL) problem is defined as follows:  Given the positions of  the anchor nodes of the WSN
and the measurable information among each non-anchor node and its neighbors,  find the positions of the rest of non-anchor nodes. 

To better describe a WSN and each sensor's possible and ideal localization actions,
game theory  is found useful in modeling  WSNs and SNL problems
\cite{jia2013distributed,bejar2010cooperative,chen2023global}. The Nash equilibrium (NE) is a prominent  concept in game theory, which characterizes a  profile of stable strategies  where rational sensor nodes would not choose to deviate from their location strategies \cite{nash1951non,belgioioso2022distributed,xu2022efficient}. Particularly, 
potential game is well-suited to model the strategic
behavior in SNL problems \cite{jia2013distributed,ke2017distributed}.
Note that the sensors need to consider the positioning accuracy of the whole WSN while ensuring their own positioning accuracy through the given information. The potential game framework can
 guarantee such an alignment between the individual sensor’s profit and the global network's objective by characterizing a global unified potential function.  In this way,  it is natural and essential to seek a global NE of the whole sensor network rather than local NE and approximate solutions, since a global NE is equal to a global optimum of the potential function denoting the network’s precise localization.

{Nevertheless, non-convexity is an intrinsic challenge of SNL problems, which cannot be avoided by selecting modeling methods.
It is the status quo that finding the global optimum or equilibrium in non-convex SNL problems is still an open problem
\cite{wang2006further,tseng2007second,wan2019sensor}.  The existing research methods for SNL problems mostly provide local or approximate solutions. 
Some relaxation methods such as  semi-definite programming  (SDP) \cite{wang2006further}  and second-order cone programming \cite{tseng2007second} are employed to transform the non-convex original %
problem into a convex optimization. 
They ignore the non-convex  constraints, yielding only approximate solutions. The alternating rank minimization (ARMA) algorithm \cite{wan2019sensor} has been considered to obtain an exact solution by mapping the rank constraints into complementary constraints. Nevertheless, this technique only guarantees the
local convergence.

In this paper, we aim to seek global solutions for SNL problems.
Specifically,
we formulate a non-convex SNL potential game, 
where 
both payoff function and potential function are characterized by continuous fourth-order polynomials. This formualtion
enables us to avoid the non-smoothness in \cite{jia2013distributed,ke2017distributed}, so as to effectively deal with the non-convex structures therein.
We reveal the existence and uniqueness of the global NE, which represents the global localization solution to  SNL. Moreover, we employ the canonical duality
	theory to transform the non-convex 
	game into a complementary dual problem and design a conjugation-based algorithm to compute the stationary points therein. Then,
 we provide a sufficient 
	condition to identify the global NE: the stationary point to the proposed algorithm 
	is the global NE if a  duality relation is satisfied. Finally, we illustrate the
effectiveness of 
our approach by numerical simulation results.

\section{Problem formulation}

In this section, we first introduce the  range-based SNL problem of interest and then formulate it as a potential game.

{
 Consider a static sensor network
	 in $\mathbb{R}^{n}$  ($ n = 2  $ or 3)  composed of 
	  $ M $ anchor nodes
whose positions are known  and $N$ non-anchor sensor nodes whose positions are
unknown (usually $ M<N $). Let a  graph $ \mathcal{G}=(\mathcal{N}, \mathcal{E}) $ represent the sensing relationships between sensors, where $\mathcal N$ is the sensor node set and $\mathcal E\subseteq\mathcal N \times \mathcal N$ is the edge set between sensors. Specifically,  
$ \mathcal{N}= \mathcal{N}_{s}\cup \mathcal{N}_{a}$, where  $ \mathcal{N}_{s}=\{1,2,\dots,N\} $ and $ \mathcal{N}_{a}= \{N+1,N+2,\dots,N+M \} $ correspond to the sets of non-anchor nodes and anchor nodes, respectively.  Let $ x_{i}^{\star} \in \mathbb{R}^{n} $ for $i\in \mathcal{N}_{s}$ denote the actual position of the $ i $-th  non-anchor
node, {and	
$ x_{N+k}^{\star} \in \mathbb{R}^{n} $ for $k\in \{1,2,\dots,M \}$
denote the actual position of anchor node $N+k\in \mathcal{N}_{a}$.}
 For a pair of sensor nodes $ i$ and $ j $, their Euclidean distance
is denoted as $ d_{ij} $. 
Each sensor has the capability
of sensing range measurements from other sensors within a
fixed range $R_s$, and $ \mathcal{E}=\{(i,j)\in\mathcal{N}\times\mathcal{N}:\|x_{i}^{\star}-x_{j}^{\star}\|\leq R_s,i\neq j\}\cup \{(i,j)\in\mathcal{N}_a\times\mathcal{N}_a,i\neq j\} $ define the edge set, i.e., there is an edge between two nodes if and only if either they are neighbors or they are both anchors.
Denote $\mathcal{N}_{s}^{i}$ as  the neighbor set of  non-anchor 
nodes $ j \in \mathcal{N}_{s}$ with $(i,j)\in\mathcal{E}$. Also, suppose that 
the measurements $d_{ij}$ are noise-free and all anchor positions $ x_{l} $, $l\in \mathcal{N}_{a}$ are accurate.
 { 	 Here we formulate the SNL  problem as an $N$-player SNL potential game $G=\{\mathcal{N}_{s}, \{\Omega_i\}_{i\in\mathcal{N}_{s}}, \{J_{i}\}_{i\in\mathcal{N}_{s}}\}$,
 where $\mathcal{N}_{s}=\{1,\dots,N\}$ corresponds to the player set,
$\Omega_i$ is player $i$'s local feasible set, which is convex and compact, and $J_{i}$ is player $i$'s payoff function.  In this context, 
we map the position estimated by  each non-anchor node as each player's strategy,
i.e.,  the strategy of the player $i$ (non-anchor node) is the estimated position $x_i\in \Omega_i$.
Denote $   \boldsymbol{\Omega}\triangleq\prod_{i=1}^{N}\Omega_{i} \subseteq \mathbb{R}^{nN} $, $ \boldsymbol{x}\triangleq \operatorname{col}\{x_{1}, \dots, x_{N}\} \in \boldsymbol{\Omega} $ as the position estimate strategy profile for all players,  and $ \boldsymbol{x}_{-i}\triangleq \operatorname{col}\{x_{1}, \dots ,x_{i-1}, x_{i+1}, \dots, x_{N}\}\subseteq \mathbb{R}^{n(N-1)}$ as the position estimate strategy profile for all players except player $ i $. 
For $i\in\mathcal{N}_s$, the payoff function $J_{i}$ is constructed as
 \begin{equation*}\label{jI}
 J_{i}(x_{i},\boldsymbol{x}_{-i})\!=\!\sum\nolimits_{j\in\mathcal{N}_{s}^{i} }(\|x_{i}-x_{j}\|^{2}\!-d_{ij}^2)^2,
 \end{equation*}
  where $(\|x_{i}\!-\!x_{j}\|^{2}\!-d_{ij}^2)^2$ in $J_{i}$  measures the   localization accuracy between  node $i$ and its   neighbor $j\!\in\!\mathcal{N}_{s}^{i}$. 

 The  individual objective of each non-anchor node  is to ensure its position accuracy,
 i.e.,
\begin{equation}\label{f1}
	\min \limits_{x_{i} \in \Omega_{i}} J_{i}\left(x_{i}, \boldsymbol{x}_{-i}\right). \quad 
\end{equation}



In the  SNL problem, each non-anchor node needs to consider the location accuracy of the whole sensor network while ensuring its own positioning accuracy through the given information. 
In other words, each non-anchor node needs to guarantee consistency between its  individual objective and   collective objective.
To this end, by regarding the individual payoff  $J_{i}$ as a marginal contribution to the whole network's collective objective \cite{marden2009cooperative,jia2013distributed}, we consider the following measurement of 
 the overall performance of sensor nodes 
\begin{equation}\label{potential-fun}
	P\!\left(x_{1},\!\dots,\! x_{N}\right)\!=\sum\nolimits_{(i,j)\in\mathcal{E} }(\|x_{i}-x_{j}\|^{2}\!-d_{ij}^2)^2.
\end{equation}
{Here,  $J_{i}$ denotes the   localization accuracy of  node $i$, which depends on the  strategies of $i$'s neighbors, while  $	P$ denotes the   localization accuracy of the entire network $\mathcal{G}$.  Then we introduce the concept of potential game.

\begin{mydef}[potential game  \cite{monderer1996potential}]\label{d2}
	A game $G=\{\mathcal{N}_{s}, \{\Omega_i\}_{i\in\mathcal{N}_{s}}, \{J_{i}\}_{i\in\mathcal{N}_{s}}\}$ is a potential
	game if there exists a potential function $P$ such that, for $i\in\mathcal{N}_s$,
	\begin{equation}\label{pp1}
	P(x_{i}^{\prime}, \boldsymbol{x}_{-i})-P\left(x_{i}, \boldsymbol{x}_{-i}\right)=J_{i}(x_{i}^{\prime}, \boldsymbol{x}_{-i})-J_{i}\left(x_{i}, \boldsymbol{x}_{-i}\right),
	\end{equation}
	for every  $\boldsymbol{x}\in \boldsymbol{\Omega}$,  and unilateral deviation $x_{i}^{\prime}\in \Omega_i$. 
\end{mydef}
It follows from Definition   \ref{d2} that 
any unilateral deviation from a strategy profile always results in the same change in both individual payoffs and a unified potential function. This indicates the alignment between each non-anchor node's selfish individual goal and the whole network's objective.

Then we verify that  $P$ in (\ref{potential-fun}) satisfies  the potential function in Definition \ref{d2}. See \cite[Appendix]{xu2024global} for the proof. 

{
\begin{proposition}\label{t3} 
With function $P$ in (\ref{potential-fun}) and payoffs $J_{i}$ for $i\in\mathcal{N}_s$ in \eqref{f1},  the   game ${{G}}=\{\mathcal{N}_{s}, \{\Omega_i\}_{i\in\mathcal{N}_{s}}, \{J_i\}_{i\in\mathcal{N}_{s}}\}$ is a  potential  game. 
%

\end{proposition}

Moreover, to attain an optimal value for $J_{i}\left(x_{i}, \boldsymbol{x}_{-i}\right)$,   players need to engage in negotiations and alter their optimal strategies. The best-known concept that describes an acceptable result achieved by all players is the NE,  whose definition is formulated below. 
\begin{mydef}[Nash equilibrium  \cite{nash1951non}]\label{d1}
	A profile $ \boldsymbol{x}^{\star}=\operatorname{col}\{x_{1}^{\star}, \dots, x_{N}^{\star} \} \in \boldsymbol{\Omega} \subseteq \mathbb{R}^{nN}$ is said to be a Nash equilibrium (NE) of game (\ref{f1}) if for any $ x_{i}\in \Omega_{i}$  we have
	\begin{equation}\label{ne}
		J_{i}\left(x_{i}^{\star}, \boldsymbol{x}_{-i}^{\star}\right) \leq J_{i}\left(x_{i}, \boldsymbol{x}_{-i}^{\star}\right),\;\forall i \in \mathcal{N}_s.
	\end{equation}
\end{mydef}
It follows from Definition \ref{d2} that
an NE of a potential game ensures not only  that each non-anchor node can adopt its optimal location strategy from the individual perspective, but also that  the sensor network as a whole can achieve a precise localization  from the global perspective.
Here, we call  NE as \textit{global} NE due to the non-convex SNL formulation in this paper. This is different from the concept of \textit{local} NE \cite{jin2020local,heusel2017gans}, which  only satisfies condition \eqref{ne} within a small neighborhood of $x_{i}^{\star}$ for $i \in \mathcal{N}_s$, rather than  the whole $\Omega_{i}$. We also  consider another mild but
well-known concept 
to help characterize the solutions to (\ref{f1}).
\begin{mydef}[Nash stationary  point \cite{pang2011nonconvex}]\label{e234}	
	A strategy profile $ \boldsymbol{x}^{\star} $ is said to be a {Nash stationary point}  of  (\ref{f1}) if
	\begin{equation}\label{e234}
	\begin{aligned}
	&\mathbf{0}_{n} \in \nabla_{x_{i}} J_{i}(x_{i}^{\star}, \boldsymbol{x}_{-i}^{\star})+\mathcal{N}_{\Omega_{i}}({x_{i}^{\star}}),  \forall i\in \mathcal{N}_s, 
	\end{aligned}
	\end{equation}
 {where $\mathcal{N}_{\Omega_{i}}(x_{i}^{\star})=\{e\in\mathbb R^n:e^T(x-x_{i}^{\star})\leq 0,\forall x\in\Omega_{i}\}$ is the normal cone at point $x_{i}^{\star}$ on set $\Omega_{i}$.}
\end{mydef}
It is not difficult to reveal that in non-convex games, if $ \boldsymbol{x}^{\star} $ is a global NE, then it must be a NE stationary point, but not vice versa.

Next, we show that global NE $\boldsymbol{x}^{\star}$ is unique and  represents the actual position profile of all non-anchor nodes, which is equal to the global solution of the SNL.
We first consider {an $n$-dimensional  representation of sensor network graph $\mathcal{G}$, which is a mapping of  $\mathcal G(\mathcal N, \mathcal E)$  to the point 
formations $\bar{\boldsymbol{x}}: \mathcal N \rightarrow \mathbb{ R}^{n}$, where 
$\bar{\boldsymbol{x}}(i)=x_i^{\operatorname{T}}$ is  the row vector of the coordinates of
the $i$-th node in $\mathbb{ R}^{n}$ and $x_i\in \mathbb{ R}^{n}$.  In this paper, the  $x_i$ is the actual position of sensor node $i$.
Given the graph  $\mathcal G(\mathcal N, \mathcal E)$ and an $n$-dimensional representation $\bar{\boldsymbol{x}}$ of it, the pair $(\mathcal G, \bar{\boldsymbol{x}})$ is called a $n$-dimensional framework. A framework  $(\mathcal G, \bar{\boldsymbol{x}})$ is called generic\footnote{Some special configurations exist among the sensor positions, e.g., groups of sensors may be collinear. The reason for using the term generic is to highlight the need to exclude the problems arising from such configurations.} if the set containing the coordinates of all its points is algebraically independent over the
rationales \cite{anderson2010formal}.  {A framework $(\mathcal G, \bar{\boldsymbol{x}})$ is called  rigid if there exists a sufficiently small positive constant $\epsilon$ such that if every framework $(\mathcal G, \bar{\boldsymbol{y}})$
satisfies  $\|x_i-y_i \|\leq\epsilon$ for $i\in\mathcal N$ and 
$\| x_i-x_j\|=\|y_i-y_j\|$ for every  pair $i,j\in \mathcal N$ connected by an edge in $\mathcal E$, then 
$\| x_i-x_j\|=\|y_i-y_j\|$ holds  for any node pair $i,j\in \mathcal N$  no matter there is an edge between them. Graph
$\mathcal G(\mathcal N, \mathcal E)$  is called generically   $n$-rigid or simply rigid (in $n$ dimensions) if any generic framework 
$(\mathcal G, \bar{\boldsymbol{x}})$ is rigid. A framework $(\mathcal G, \bar{\boldsymbol{x}})$  is globally rigid if every framework $(\mathcal G, \bar{\boldsymbol{y}})$ satisfying
$\| x_i-x_j\|=\|y_i-y_j\|$  for any node pair $i,j\in \mathcal N$  connected by an edge in $\mathcal E$ and  $\| x_i-x_j\|=\|y_i-y_j\|$  for any node pair $i,j\in \mathcal N$  that are not connected by a single edge. 
Graph
$\mathcal G(\mathcal N, \mathcal E)$  is called generically globally rigid if 
any generic framework 
$(\mathcal G, \bar{\boldsymbol{x}})$  is 
globally rigid \cite{anderson2010formal,fidan2010closing,tay1985generating}.}
On this basis, we make the following basic assumption. 

\begin{assumption}
 The sensor topology graph	 $ \mathcal{G} $ is undirected and generically globally rigid.
\end{assumption}
The undirected graph topology is usually a common assumption in many graph-based approaches \cite{chen2021distributed,jing2021angle}. 
The connectivity of $\mathcal{G} $ can also be induced by some disk graph \cite{wan2019sensor},  which ensures the validity of the information transmission between nodes.
{The generic global rigidity of $\mathcal{G} $ has been widely employed in SNL problems to guarantee the graph structure invariant, which indicates a unique localization of the sensor network  \cite{calafiore2010distributed,shi2010distributed,anderson2008rigid}.
Besides, there have been extensive discussions on graph rigidity in existing works \cite{wan2019sensor,cao2021bearing}, but it is not the primary focus of our paper. }

The following lemma reveals the existence and uniqueness of global NE $\boldsymbol{x}^{\star}$. See \cite[Appendix]{xu2024global} for the proof. 

{
\begin{lemma}\label{l11} 
	Under Assumption 1, the global NE $\boldsymbol{x}^{\star}$ of the potential game G is unique and corresponds to the actual position profile of all non-anchor nodes, which represents the global solution of the SNL.
\end{lemma}


 
%
%
 
While we have obtained guarantees regarding the existence and uniqueness of global NE of the SNL problem, its identification and computation are still challenging since $J_{i}$ and $	P$  are non-convex functions in our model. 
Actually, as for convex games, most of the existing research works seek global NE   via investigating first-order stationary points  under Definition 3 \cite{facchinei2010penalty,koshal2016distributed,chen2021distributed}.
 However, in such a non-convex regime (\ref{potential-fun}), 
one cannot expect to find a global NE easily following this way, because  stationary points in non-convex settings are not equivalent to  global NE anymore. 
Such similar potential game models  have  also been considered in \cite{jia2013distributed,ke2017distributed}.
As different from the use of the Euclidean norm in \cite{jia2013distributed,ke2017distributed}, i.e., $\|\|x_{i}-x_{j}\|\!-d_{ij}\|$, 
we adopt the square of Euclidean norm to characterize  $J_{i}$ and $P$, i.e., $\|\|x_{i}-x_{j}\|^2-d_{ij}^2\|^2$. These functions endowed with continuous fourth-order polynomials enable us to avoid the non-smoothness and deal with the inherent non-convexity of SNL with useful technologies, so as to get the global NE.  
On the other hand, previous efforts merely yield an approximate solution or a local NE by relaxing non-convex constraints or relying on additional convex assumptions,  either under potential games 
or other modeling methods \cite{calafiore2010distributed,wan2019sensor}. 
Thus, they
fail to adequately address the intrinsic 
non-convexity of SNL.}
To this end, we investigate the
identification condition of the global NE in the SNL problem. Specifically, we aim to find
the conditions that a stationary point of (\ref{f1}) is consistent with the global NE and design an algorithm to solve it.

\section{Derivation of the global Nash equilibrium}\label{sec:duality}

{
In this section, we   explore 
the identification condition of the global NE of the SNL problem by virtue of canonical dual theory and develop a conjugation-based algorithm to compute it. }  



%


  
 It is hard to directly identify  whether a stationary point is the global NE
  on  the non-convex potential function (\ref{potential-fun}).  Here, 
we employ canonical duality theory \cite{gao2017canonical} to transform (\ref{potential-fun}) into a complementary dual problem  and  investigate the relationship between  a stationary point of the dual problem and the global NE of game (\ref{f1}). 

\noindent\textbf{Canonical transformation}
We first reformulate  (\ref{potential-fun}) in a canonical form.
Define $\xi_{ij}={\Lambda_{ij}(\boldsymbol{x})}=\|x_{i}-x_{j}\|^2$.
 in 
 (\ref{potential-fun})
  and define the profiles 
  \begin{equation}\label{ci} 
  		{\Lambda}(\boldsymbol{x})=\operatorname{col}\{\Lambda_{ij}(\boldsymbol{x})\}_{(i,j)\in\mathcal{E}},
  \;{\xi}=\operatorname{col}\{\xi_{ij}\}_{(i,j)\in\mathcal{E} }\in \Xi.
\end{equation}
 Here,
 ${\Lambda}(\boldsymbol{x})$ 
 map the decision variables in domain $\boldsymbol{{\Omega}}$ to the quadratic functions in space $\Xi \subseteq\mathbb{R}^{|\mathcal{E}| } $.
Moreover, we introduce quadratic functions $\Phi:  \Xi\rightarrow \mathbb{R}$,
\begin{equation}\label{cl}
	\Phi({\xi})\!= \sum\nolimits_{(i, j) \in \mathcal{E}}(\xi_{i j}-\!d_{i j}^2)^2.
\end{equation}
Thus, the potential function (\ref{potential-fun}) can  be  rewritten as:
$
P(\boldsymbol{x})=\Phi({\Lambda}(\boldsymbol{x})). 
$
{ Note that the gradients $\nabla\Phi: \Xi\rightarrow \Xi^{*}$ is  a one-to-one mapping, where  $\Xi^{*}$ is the range space of the  gradient. Thus, recalling \cite{gao2017canonical}, 
	 $ \Phi: \Xi \rightarrow \mathbb{R}$ is a  convex differential
canonical function. }
This indicates that
 the following one-to-one duality relation is invertible on $\Xi \times \Xi^*$:
\begin{equation}\label{dual}
	\begin{array}{ll}
		\tau_{i j}=\nabla_{\xi_{i j}}\Phi({\xi})=2(\xi_{i j}-d_{i j}^2), & (i, j) \in \mathcal{E}.
	\end{array}
\end{equation}
 Denote the profiles $\tau=\operatorname{col}\{\tau_{ij}\}_{(i,j)\in\mathcal{E} }\in\Xi^{*}\subseteq \mathbb{R}^{q}$,
 where {$q=\left|\mathcal{E}\right|$ is the total number of elements in the edge sets $\mathcal{E}$. 
Based on \eqref{dual}, the Legendre conjugates of $\Phi$ 
can be uniquely defined by 
 \begin{equation}\label{cr}
 	\begin{aligned}
 		\Phi^*({\tau}) &\! =\! (\xi)^{\operatorname{T}} \tau\!-\!\Phi(\xi)\!=\!
 		\sum\nolimits_{(i, j) \in \mathcal{E}} \frac{1}{4}(\tau_{i j})^2\!+\!d_{i j}^2 \tau_{i j}, 
 	\end{aligned}
 \end{equation}
 where $(\xi,\tau)$ is called the Legendre canonical duality pair
on $\Xi\times\Xi^{*}$.
We regard  $\tau$ as a canonical dual variable on the  dual space $\Xi^{*}$. Then, based on the canonical duality theory  \cite{gao2017canonical}, we define the following the complementary function  $\Psi: \boldsymbol{\Omega}\times \Xi^{*}\rightarrow \mathbb{R} $,
	\begin{align}\label{complementary}
		\Psi\left(x_{1},\!\dots,\! x_{N},\tau\right)=&(\xi)^{\operatorname{T}} \tau-\Phi^*(\tau)\\
		=&\!\!\sum_{(i,j)\in\mathcal{E} }\!\!\!\!\tau_{ij}(\|x_{i}-x_{j}\|^{2}-d_{ij}^2)\!-\!\!\sum_{(i,j)\in\mathcal{E}}\!\!\!\!\frac{(\tau_{ij})^2}{4}. \nonumber
	\end{align}
So far, we have transformed the non-convex function (\ref{potential-fun}) into a complementary dual problem (\ref{complementary}).
We have the following result about the equivalency relationship of stationary points between (\ref{complementary}) and  (\ref{potential-fun}), whose proof is shown in Appendix \ref{bb}.
\begin{mythm}\label{t1}
 { For a profile $ \boldsymbol{x}^{\star} $, if  there exists ${\tau}^{\star}\in {\Xi}^{*}$  such that  for $ i\in \mathcal{N}_s$, 
	$ (\boldsymbol{x}^{\star},{\tau^{\star}})$
	is a stationary point of  complementary function  $ \Psi(\boldsymbol{x},\tau)$, then $ \boldsymbol{x}^{\star} $ is a Nash stationary point of game (\ref{f1}). }
\end{mythm}

By Theorem \ref{t1}, the equivalency of   stationary points between (\ref{complementary}) and  (\ref{f1})
is due to the fact that
 the duality relations (\ref{dual}) are  unique and   invertible on $ \boldsymbol{\Omega}\times\Xi^{*} $,   thereby closing the duality gap between  the non-convex original game and its canonical dual problem.

\noindent\textbf{Sufficient feasible domain}
 Next, we introduce a sufficient feasible domain for the introduced conjugate variable $\tau$, in order to investigate the global optimality of the stationary points in (\ref{complementary}).
Consider the  second-order derivative
of $\Psi\left(\boldsymbol{x},\tau\right)$ in $\boldsymbol{x}$. Due to the expression of (\ref{complementary}), we can find that $\Psi$ is quadratic in $\boldsymbol{x}$. Thus, $\nabla^{2}_{\boldsymbol{x}}\Psi$ is   $\boldsymbol{x}$-free,  and is indeed a linear combination for the elements of $\tau$.
In this view, we denote $	Q(\tau)= \nabla^{2}_{\boldsymbol{x}}\Psi$.   On this basis, we introduce the following set of $\tau$ 
\begin{equation}\label{s23}
	\mathbb{E}^{+}= \Xi^{*}\cap \{\tau : 
	Q(\tau) \succeq \boldsymbol{0}_{nN}
	\}. 
\end{equation} 
}

\noindent\textbf{Algorithm design}\label{sec:cent}}
Then, we design a  conjugation-based algorithm to compute the stationary points of the SNL problem with the assisted complementary information (the   Legendre conjugate of $\Phi$  and the canonical conjugate variable $\tau$).

In Alg. 1,  the terms about $-\nabla_{x_{i}} \Psi( \boldsymbol{x}[k],\tau[k])$ for $i\in\mathcal{N}_s$ and $\nabla_{\tau} \Psi( \boldsymbol{x}[k],\tau[k])$ represent the directions of gradient descent and ascent according to $\Psi$. The terms about $\Pi_{\mathbb{E}^{+}}$ and $\Pi_{\Omega_{i}}$ are projection operators \cite{xu2023algorithm}.
When $ \tau\in\mathbb{E}^{+}  $,  the positive semi-definiteness of $Q(\tau)$ implies that $\Psi(\boldsymbol{x},{\tau}  )$ is convex with respect to $\boldsymbol{x}$. Besides, the convexity of $\Phi(\xi)$ derives that its Legendre conjugate $\Phi^{*}(\tau)$ is also convex \cite{10044227}, implying that the  complementary function  $\Psi(\boldsymbol{x},{\tau}  )$  is concave in $\tau$. Together with the non-expansiveness of projection operators and a decaying step size $\{\alpha[k]\}$, this convex-concave property of $\Psi$ implies the convergence of Alg. 1 and enables us to  identify the global NE.

\begin{algorithm}[t]
	\caption{Conjugation-based SNL algorithm}
	\label{centralized-alg}
	\renewcommand{\algorithmicrequire}{\textbf{Input:}}
	\renewcommand{\algorithmicensure}{\textbf{Initialize:}}
	\begin{algorithmic}[1]
		\REQUIRE Step size  $ \{\alpha[k] \} $.
		\ENSURE Set $\tau[0]\in {\mathbb{E}^{+}}, x_{i}[0]\in \Omega_{i}, \,i\in\{1,\dots,N\}  $, 
		\FOR{$k = 1,2,\dots$}
		\STATE  update the shared  canonical dual variable:\;\vspace{0.1cm}
		
		$\tau[k+1] = \Pi_{\mathbb{E}^{+}}(\tau[k] +\alpha[k]
		\nabla_{\tau} \Psi( \boldsymbol{x}[k],\tau[k])	
		)
		$
		\FOR{$i = 1,\dots N$}
		\STATE update the decision  variable of non-anchor node $i$:\;\vspace{0.1cm}
		
		$	x_{i}[k+1] =\Pi_{\Omega_{i}}( x_i[k]-\alpha[k] \nabla_{x_{i}} \Psi( \boldsymbol{x}[k],\tau[k]) )
		$
		\ENDFOR
		\ENDFOR
	\end{algorithmic}
\end{algorithm}
\noindent\textbf{Equilibrium design}\label{sec:eq}}
On this basis,
 we  establish the relationship between the global NE in  (\ref{potential-fun})
and a stationary point computed from Alg. 1. The proof is shown in Appendix \ref{cc}.

{
	\begin{mythm}\label{t2}
	Under Assumption 1,	
	profile $\boldsymbol{x}^{\star}$ is the global NE		of game ${{G}}$ if  there exists ${\tau}^{\star}\in \mathbb{E}^{+}$ such that a  stationary point 
	$(\boldsymbol{x}^{\star},{\tau}^{\star}) \in \boldsymbol{\Omega}\times\mathbb{E}^{+} $  obtained from Alg. 1 
	  satisfies
	$$\tau_{i j}^{\star}=\nabla_{\xi_{i j}}\Phi({\xi})|_{\xi_{i j}=\|x_{i}^{\star}-x_{j}^{\star}\|^2}, \forall (i, j) \in \mathcal{E}. $$  
\end{mythm}}

The result in Theorem \ref{t2}  reveals that  
once the stationary point of Alg. 1 is obtained, we can check the duality  relation $\tau_{i j}^{\star}=\nabla_{\xi_{i j}}\Phi({\xi})|_{\xi_{i j}=\|x_{i}^{\star}-x_{j}^{\star}\|^2}, \forall (i, j) \in \mathcal{E} $,
 so as to identify whether the solution of Alg. 1  is the global NE.
In fact, 	it is necessary to check the duality for the convergent point $(\boldsymbol{x}^{\star},\tau^{\star})$ of Alg. 1, because the computation of   $\tau^{\star}$ is restricted on the sufficient domain $\mathbb{E}^{+}$ instead of the original $\Xi^*$. In this view, 	
	the gradient of $\tau^{\star}$ may fall into the normal cone $\mathcal{N}_{\mathbb{E}^{+}}(\tau^{\star})$ instead of being equal to $\boldsymbol{0}_{q}$, thereby losing the one-to-one relationship with $\boldsymbol{x}^{\star}$. Thus,  $\boldsymbol{x}^{\star}$ may not be the global NE.
	In addition, we cannot directly employ the standard Lagrange multiplier
	method and the associated Karush-Kuhn-Tucker (KKT) theory herein, because we need to first confirm a feasible domain of $\tau$ by utilizing
	canonical duality information (referring to $\Xi^*$). In other words, once the duality relation is verified, we can say that the convergent point $(\boldsymbol{x}^{\star},\tau^{\star})$ of Alg. 1 is indeed the global NE of game \eqref{f1}.

We summarize a road map for seeking global
	NE in this non-convex SNL problem for friendly comprehension. That is, 
	once
	the problem is defined and formulated,  we first transform the original SNL potential game into a dual complementary problem. Then we seek the stationary point of  $ \Psi (\boldsymbol{x},\tau)$  via algorithm iterations, wherein the dual variable $\tau$
	is restricted on  $\mathbb{E}^{+}$. Finally, after
	obtaining the stationary point by convergence, we
 identify whether the convergent point satisfies the
	duality relation. If so, the
	convergent point is the global NE.
\section{Numerical Experiments}

In this section, 
we examine the effectiveness of our approach to seek the
global NE of the SNL problem.  


    \begin{figure*}[bp]
			\hspace{-0.7cm}
			\centering	
			\subfigure[$N=10$]{
				\begin{minipage}[t]{0.248\linewidth}
					\centering
					\includegraphics[width=4.7cm]{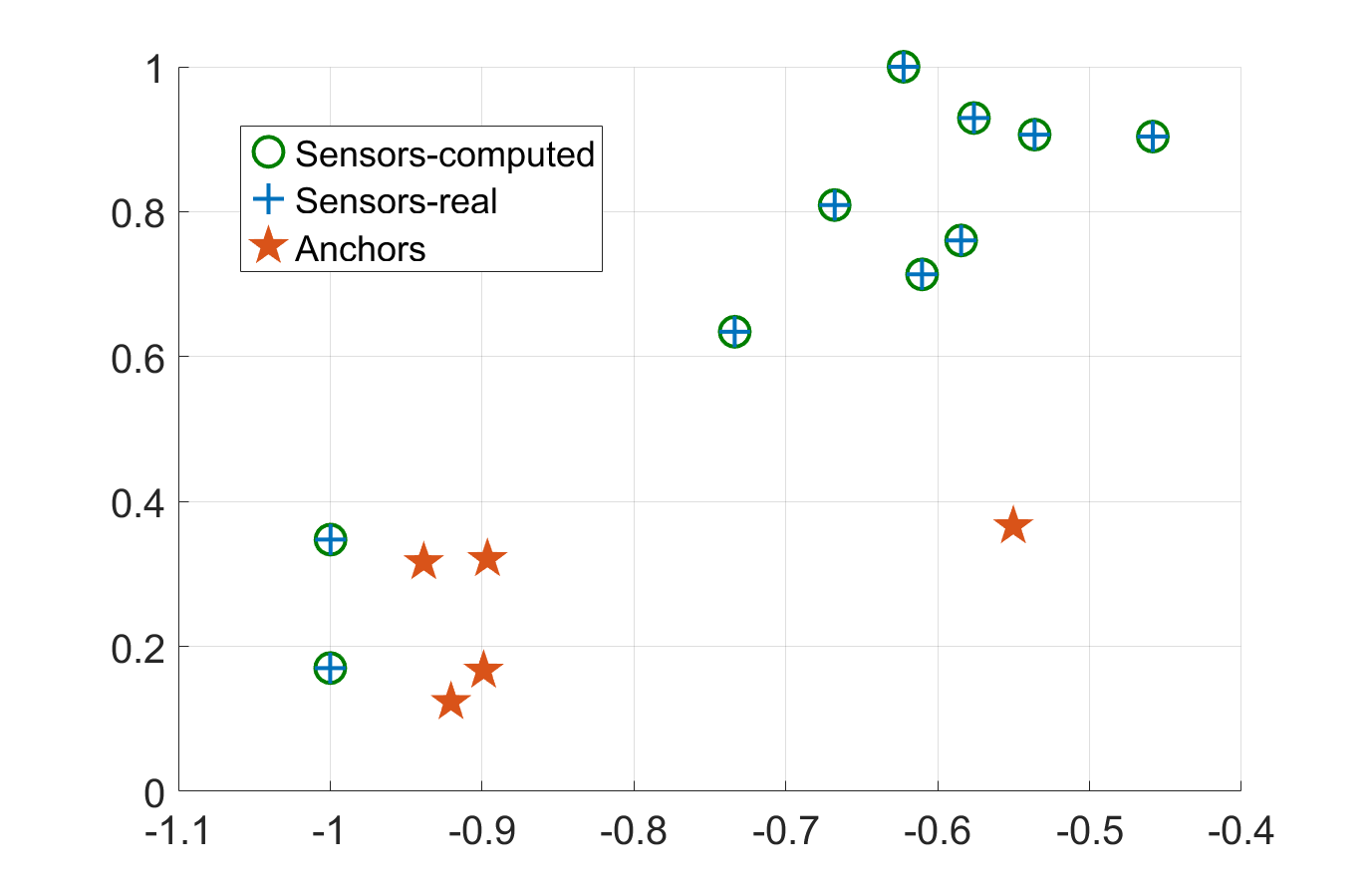}
				\end{minipage}%
			}%
			\subfigure[$N=20$]{
				\begin{minipage}[t]{0.248\linewidth}
					\centering
					\includegraphics[width=4.7cm]{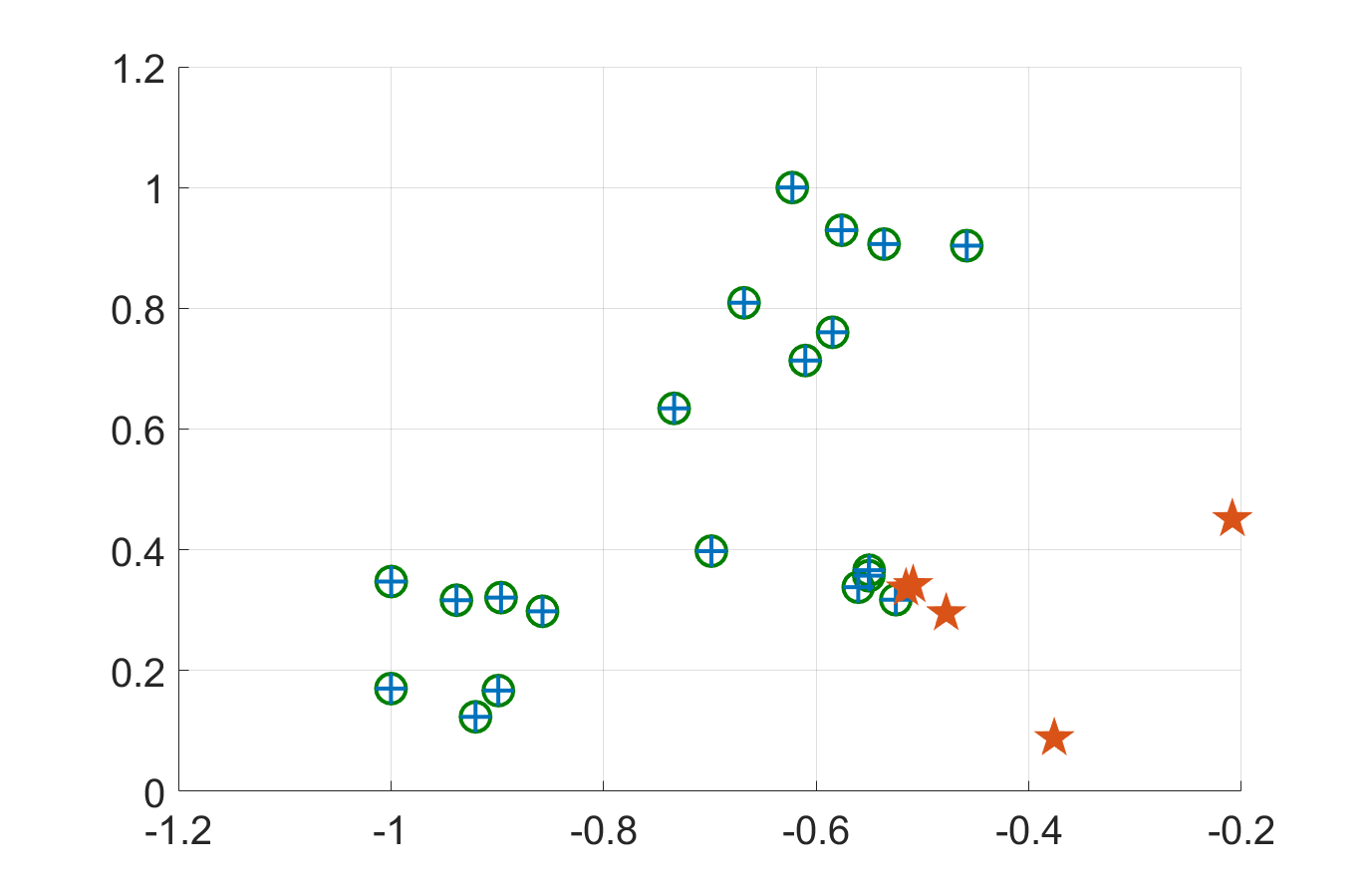}
				\end{minipage}%
			}%
   \subfigure[$N=35$]{
				\begin{minipage}[t]{0.248\linewidth}
					\centering
					\includegraphics[width=4.7cm]{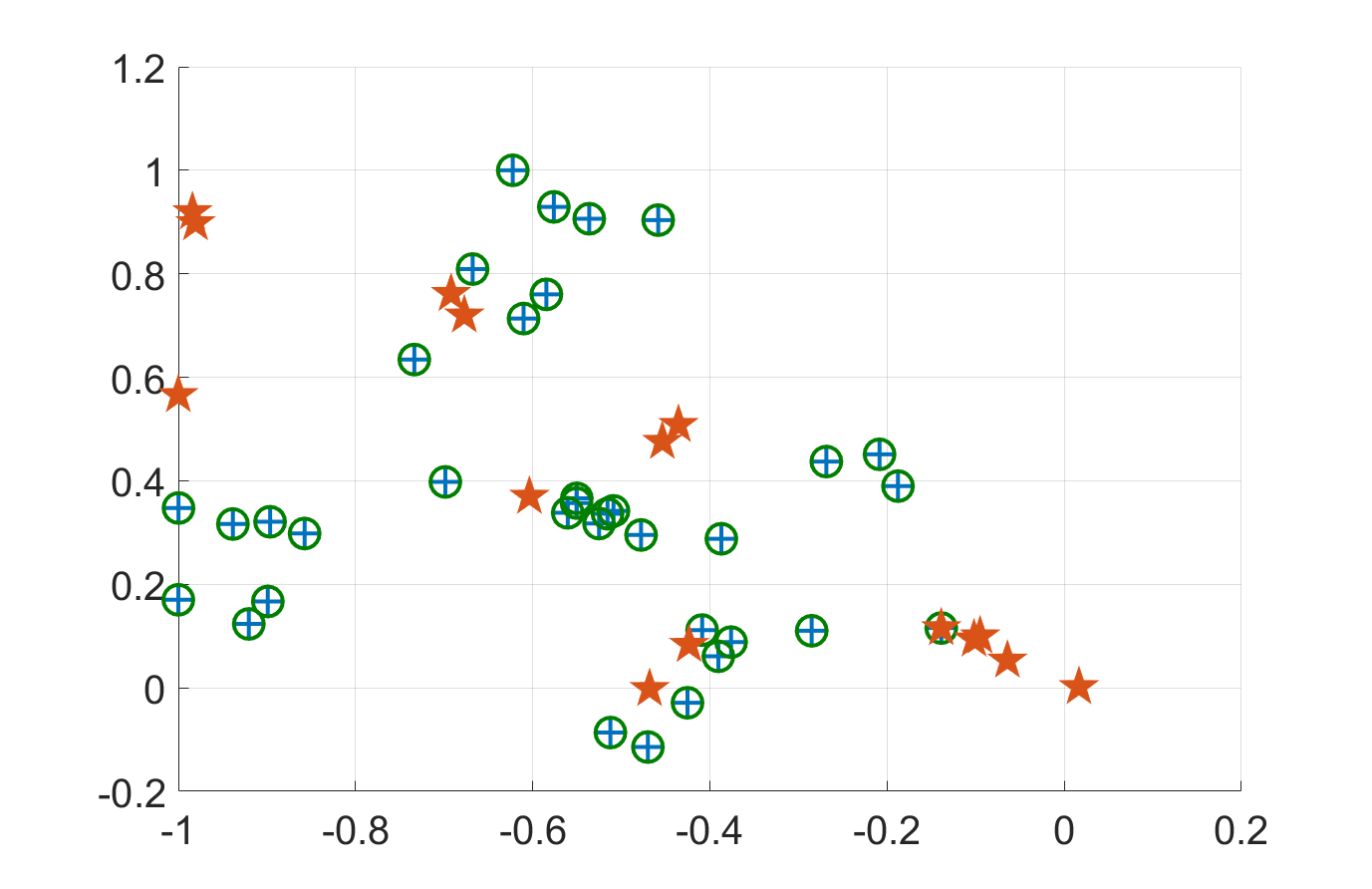}
     \vspace{-0.2cm}
				\end{minipage}%
			}%
			\subfigure[$N=50$]{
				\begin{minipage}[t]{0.248\linewidth}
					\centering
					\includegraphics[width=4.7cm]{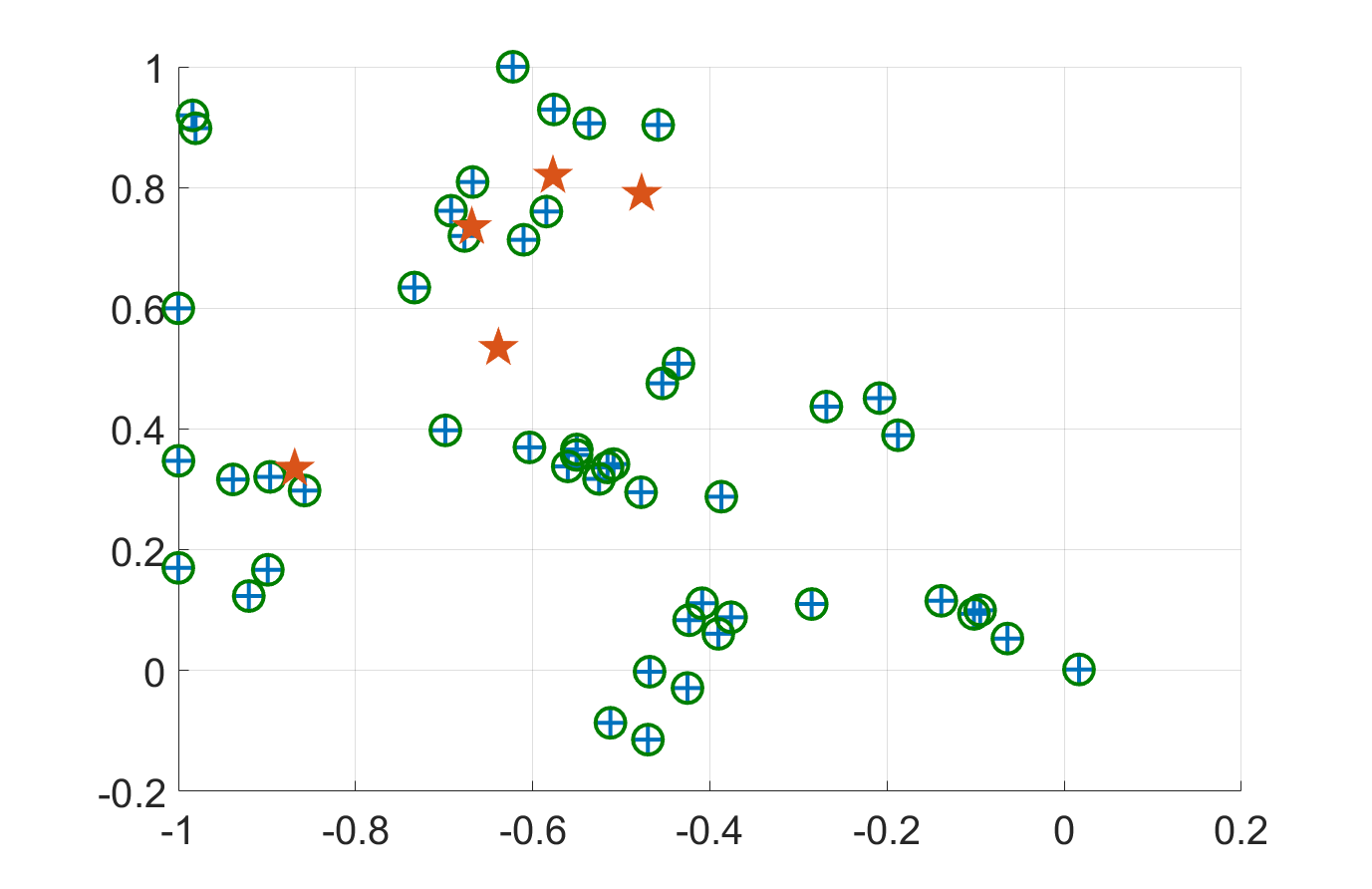}
				\end{minipage}%
			}
			\centering
			\caption{Computed sensor location results with different configurations.} 
			\label{figg4}
		\end{figure*}
We first  consider a two-dimensional case based on the UJIIndoorLoc dataset. The UJIIndoorLoc dataset was introduced in 2014 at the International Conference on Indoor Positioning and Indoor Navigation, to estimate a user location based on building and floor. The dataset is available on the UC Irvine Machine Learning Repository website \cite{23}. We extract the latitude and longitude coordinates of part of the sensors and standardize the data by doing min-max normalization.  
 We employ Alg. 1 to solve this problem. 
Set the  tolerance $t_{tol}= 10^{-5}$  and the terminal criterion 
$ \|\boldsymbol{x}[k+1]-\boldsymbol{x}[k]\| \leq t_{tol},\; \|{\tau}[k+1]-{\tau}[k+1]\| \leq t_{tol}.  $

We show the effectiveness of Alg. 1 for SNL problems with different node configurations.  
Take  $ N=10,20,35,50 $ and different numbers of anchor nodes. 
Fig. \ref{figg4} shows the computed sensor location results in these cases. The anchor nodes and
the true locations of non-anchor nodes are  shown by red
stars and blue asterisks, and the computed locations are shown
by green circles.
We can see that  Alg. 1  can  localize all sensors in either small or large sensor network sizes.

\section{Conclusion}
In this paper, we have focused on the non-convex SNL problems. We have presented novel results on the identification condition of the global solution and the position-seeking algorithms. By formulating a non-convex SNL potential game, we have shown that the global NE   exists and is unique. Then based on the canonical duality theory, we have proposed  a conjugation-based algorithm  to compute the stationary point of a complementary dual problem, which  actually induces the global
NE  if a duality relation can be checked. Finally, the computational efficiency
 of our algorithm has been
illustrated by several experiments. 

In the future, we may extend our current results to more complicated cases  such as i) generalizing the algorithm to distributed situations,  ii) generalizing the model to cases
with measurement noise, and iii) exploring milder graph conditions.
\appendices
\renewcommand\thesection{\Alph{section}}

\section{Proof of Theorem  \ref{t1}}\label{bb}

	If 
	there exists ${\tau}^{\star}\in \Xi^{*}$ such that $(\boldsymbol{x}^{\star},\tau^{\star})$ is a stationary point of $ \Psi(\boldsymbol{x},\tau)$, 
	then it satisfies the first-order condition, that is 
\begin{subequations}\label{ddd}
\begin{align}
	&\mathbf{0}_{nN} \in   \nabla_{\boldsymbol{x}}\,\tau^{\star\operatorname{T}} \Lambda\left(\boldsymbol{x}^{\star}\right)	
	+\mathcal{N}_{\boldsymbol{\Omega}}(\boldsymbol{x}^{\star}),\\
	&\mathbf{0}_{q} \in -\Lambda(\boldsymbol{x}^{\star})+\nabla \Phi^{*}(\tau^{\star})+ \mathcal{N}_{\Xi^{*}}(\tau^{\star}),
 \end{align}
\end{subequations} 
Moreover, 
 based on the  invertible one-to-one duality relation \eqref{dual},  for given $\xi \in \Xi$ 
with $\xi=\Lambda(\boldsymbol{x}^{\star})$,
 we have $$\tau_{i j}^{\star}=\nabla_{\xi_{i j}}\Phi({\xi})|_{\xi_{i j}=\|x_{i}^{\star}-x_{j}^{\star}\|^2} \Leftrightarrow \xi_{ij}=\nabla \Phi^{*}(\tau_{ij}^{\star})$$ for $ (i, j) \in \mathcal{E}$. By employing this relation in (\ref{ddd}b), we have $
\mathbf{0}_{q} = -\Lambda(\boldsymbol{x}^{\star})+\nabla \Phi^{*}(\tau^{\star}),
 $
 which implies 
$
 \tau^{\star}=\nabla \Phi(\Lambda(\boldsymbol{x}^{\star})).
$
 By  substituting $\tau^{\star}$ with $\nabla \Phi(\Lambda(\boldsymbol{x}^{\star}))$, we have \begin{equation}\label{chain}
    \mathbf{0}_{nN} \in   \nabla_{\boldsymbol{x}}\, \Phi(
	\Lambda(\boldsymbol{x}^{\star}))^{\operatorname{T}} \Lambda\left(\boldsymbol{x}^{\star}\right)	
	+\mathcal{N}_{\boldsymbol{\Omega}}(\boldsymbol{x}^{\star}).
\end{equation}
According to the chain rule, 
 $
\nabla \Phi(
	\Lambda(\boldsymbol{x}^{\star}))^{\operatorname{T}} \Lambda\left(\boldsymbol{x}^{\star}\right)	= \nabla_{\boldsymbol{x}}\Phi(\boldsymbol{x}^{\star}).
 $
 Therefore, \eqref{chain} is equivalent to 
\begin{equation}\label{ft41}
	 \mathbf{0}_{nN} \in   \nabla_{\boldsymbol{x}}\Phi(\boldsymbol{x}^{\star})+\mathcal{N}_{\boldsymbol{\Omega}}(\boldsymbol{x}^{\star}).
\end{equation} 
According to the definition of potential game, \eqref{ft41} implies 
\begin{equation}\label{dde}
	\mathbf{0}_{n} \in \nabla_{x_{i}} J_{i}(x_{i}^{\star}, \boldsymbol{x}_{-i}^{\star})+\mathcal{N}_{\Omega_{i}}({x_{i}^{\star}}),
\end{equation}
 which yields the conclusion. \hfill $\square$ 

\section{Proof of Theorem  \ref{t2}}\label{cc}
{
	If 
	there exists ${\tau}^{\star}\in \mathbb{E}^{+}$ such that the pair $(\boldsymbol{x}^{\star},\tau^{\star})$ is a stationary point of Alg. 1, 
	then it satisfies the first-order condition with respect to  $ \Psi(x_{i},\boldsymbol{x}_{-i},\tau)$, that is 
\begin{subequations}\label{dds}
\begin{align}
	&\mathbf{0}_{nN} \in   \nabla_{\boldsymbol{x}}\,\tau^{\star\operatorname{T}} \Lambda\left(\boldsymbol{x}^{\star}\right)	
	+\mathcal{N}_{\boldsymbol{\Omega}}(\boldsymbol{x}^{\star}),\\
	&\mathbf{0}_{q} \in -\Lambda(\boldsymbol{x}^{\star})+\nabla \Phi^{*}(\tau^{\star})+ \mathcal{N}_{\mathbb{E}^+}(\tau^{\star}),
 \end{align}
\end{subequations}  
Together with $\tau_{i j}^{\star}=\nabla_{\xi_{i j}}\Phi({\xi})|_{\xi_{i j}=\|x_{i}^{\star}-x_{j}^{\star}\|^2}, \forall (i, j) \in \mathcal{E} $, we claim that  
	the canonical duality relation holds over $\boldsymbol{\Omega}\times\mathbb{E}^{+}$. Thus, (\ref{dds}b) becomes $
\mathbf{0}_{q} = -\Lambda(\boldsymbol{x}^{\star})+\nabla \Phi^{*}(\tau^{\star}).
 $
This indicates that 
the stationary point  $ (\boldsymbol{x}^{\star},{\tau^{\star}})$
	of    $ \Psi(x_{i},\boldsymbol{x}_{-i},\tau)$ on $\boldsymbol{\Omega}\times\mathbb{E}^{+}$ is also a stationary point profile of $\Psi$ on $\boldsymbol{\Omega}\times\Xi^{*}$. Based on Theorem \ref{t1}, we can further derive that the profile $\boldsymbol{x}^{\star}$ with respect to the stationary point  $ (\boldsymbol{x}^{\star},{\tau^{\star}})$ of $\Psi$ on $\boldsymbol{\Omega}\times\mathbb{E}^{+}$ is a Nash stationary point of game \eqref{f1}.

	Moreover, recall $
	\mathbb{E}^{+}= \Xi^{*}\cap \{\tau : 
	Q(\tau) \succeq \boldsymbol{0}_{nN}\}
 $
 with  $	Q(\tau)= \nabla^{2}_{\boldsymbol{x}}\Psi$.
 This indicates that 
 $ \Psi(\boldsymbol{x},\tau)$  is convex in $\boldsymbol{x}$. Also, note that $ \Psi(\boldsymbol{x},\tau)$ is  concave in  dual  variable $\tau$ due to the convexity of $\Phi(\cdot)$.  
 
 Thus,
	 we can obtain the  global optimality of $(\boldsymbol{x}^{\star}, \boldsymbol{\tau}^{\star})$ on $\boldsymbol{\Omega}\times \mathbb{E}^{+} $, that is, for $\boldsymbol{x}\in \boldsymbol{\Omega}$ and $\tau\in \mathbb{E}^{+}$, 
	\begin{equation*}
	\Psi(\boldsymbol{x}^{\star},  \tau)\leq\Psi(\boldsymbol{x}^{\star},  \tau^{\star}) \leq \Psi(\boldsymbol{x},  \tau^{\star}). 
	\end{equation*}
	The inequality relation above tells that 
\begin{equation*}
	J_{i}(x_{i}^{\star}, \boldsymbol{x}_{-i}^{\star})\leq J_{i}(x_{i}, \boldsymbol{x}_{-i}^{\star}), \quad \forall x_{i}\in \Omega_i,\quad \forall i\in\mathcal{N}_s.
\end{equation*}
This confirms that $\boldsymbol{x}^{\star}$ is the global NE of (\ref{f1}), which completes the proof. \hfill $\square$

	\bibliographystyle{IEEEtran}
\bibliography{reference}

\begin{thebibliography}{10}
\providecommand{\url}[1]{#1}
\csname url@samestyle\endcsname
\providecommand{\newblock}{\relax}
\providecommand{\bibinfo}[2]{#2}
\providecommand{\BIBentrySTDinterwordspacing}{\spaceskip=0pt\relax}
\providecommand{\BIBentryALTinterwordstretchfactor}{4}
\providecommand{\BIBentryALTinterwordspacing}{\spaceskip=\fontdimen2\font plus
\BIBentryALTinterwordstretchfactor\fontdimen3\font minus
  \fontdimen4\font\relax}
\providecommand{\BIBforeignlanguage}[2]{{%
\expandafter\ifx\csname l@#1\endcsname\relax
\typeout{** WARNING: IEEEtran.bst: No hyphenation pattern has been}%
\typeout{** loaded for the language `#1'. Using the pattern for}%
\typeout{** the default language instead.}%
\else
\language=\csname l@#1\endcsname
\fi
#2}}
\providecommand{\BIBdecl}{\relax}
\BIBdecl

\bibitem{akyildiz2002wireless}
I.~F. Akyildiz, W.~Su, Y.~Sankarasubramaniam, and E.~Cayirci, ``Wireless sensor
  networks: a survey,'' \emph{Computer {N}etworks}, vol.~38, no.~4, pp.
  393--422, 2002.

\bibitem{liu2015event}
Q.~Liu, Z.~Wang, X.~He, and D.-H. Zhou, ``Event-based recursive distributed
  filtering over wireless sensor networks,'' \emph{IEEE Transactions on
  Automatic Control}, vol.~60, no.~9, pp. 2470--2475, 2015.

\bibitem{marechal2010joint}
N.~Marechal, J.-M. Gorce, and J.-B. Pierrot, ``Joint estimation and gossip
  averaging for sensor network applications,'' \emph{IEEE Transactions on
  Automatic Control}, vol.~55, no.~5, pp. 1208--1213, 2010.

\bibitem{jing2021angle}
G.~Jing, C.~Wan, and R.~Dai, ``Angle-based sensor network localization,''
  \emph{IEEE Transactions on Automatic Control}, vol.~67, no.~2, pp. 840--855,
  2021.

\bibitem{sun2010corrosion}
G.~Sun, G.~Qiao, and B.~Xu, ``Corrosion monitoring sensor networks with energy
  harvesting,'' \emph{IEEE Sensors Journal}, vol.~11, no.~6, pp. 1476--1477,
  2010.

\bibitem{sun2005reliable}
T.~Sun, L.-J. Chen, C.-C. Han, and M.~Gerla, ``Reliable sensor networks for
  planet exploration,'' in \emph{Proceedings. 2005 IEEE Networking, Sensing and
  Control, 2005.}\hskip 1em plus 0.5em minus 0.4em\relax IEEE, 2005, pp.
  816--821.

\bibitem{jing2018weak}
G.~Jing, G.~Zhang, H.~W. Joseph~Lee, and L.~Wang, ``Weak rigidity theory and
  its application to formation stabilization,'' \emph{SIAM Journal on Control
  and optimization}, vol.~56, no.~3, pp. 2248--2273, 2018.

\bibitem{estrin2000embedding}
D.~Estrin, R.~Govindan, and J.~Heidemann, ``Embedding the internet:
  introduction,'' \emph{Communications of the ACM}, vol.~43, no.~5, pp. 38--41,
  2000.

\bibitem{fang2023distributed}
X.~Fang, L.~Xie, and X.~Li, ``Distributed localization in dynamic networks via
  complex {L}aplacian,'' \emph{Automatica}, vol. 151, p. 110915, 2023.

\bibitem{mao2007wireless}
G.~Mao, B.~Fidan, and B.~D. Anderson, ``Wireless sensor network localization
  techniques,'' \emph{Computer {N}etworks}, vol.~51, no.~10, pp. 2529--2553,
  2007.

\bibitem{wan2019sensor}
C.~Wan, G.~Jing, S.~You, and R.~Dai, ``Sensor network localization via
  alternating rank minimization algorithms,'' \emph{IEEE Transactions on
  Control of Network Systems}, vol.~7, no.~2, pp. 1040--1051, 2019.

\bibitem{ahmadi2021distributed}
S.~P. Ahmadi, A.~Hansson, and S.~K. Pakazad, ``Distributed localization using
  {L}evenberg-{M}arquardt algorithm,'' \emph{EURASIP Journal on Advances in
  Signal Processing}, vol. 2021, no.~1, pp. 1--26, 2021.

\bibitem{jia2013distributed}
J.~Jia, G.~Zhang, X.~Wang, and J.~Chen, ``On distributed localization for road
  sensor networks: A game theoretic approach,'' \emph{Mathematical Problems in
  Engineering}, vol. 2013, 2013.

\bibitem{bejar2010cooperative}
B.~Bejar, P.~Belanovic, and S.~Zazo, ``Cooperative localisation in wireless
  sensor networks using coalitional game theory,'' in \emph{2010 18th European
  Signal Processing Conference}.\hskip 1em plus 0.5em minus 0.4em\relax IEEE,
  2010, pp. 1459--1463.

\bibitem{chen2023global}
G.~Chen, G.~Xu, F.~He, Y.~Hong, L.~Rutkowski, and D.~Tao, ``Global {N}ash
  equilibrium in non-convex multi-player game: Theory and algorithms,''
  \emph{arXiv preprint arXiv:2301.08015}, 2023.

\bibitem{nash1951non}
J.~Nash, ``Non-cooperative games,'' \emph{Annals of Mathematics}, pp. 286--295,
  1951.

\bibitem{belgioioso2022distributed}
G.~Belgioioso, P.~Yi, S.~Grammatico, and L.~Pavel, ``Distributed generalized
  {N}ash equilibrium seeking: An operator-theoretic perspective,'' \emph{IEEE
  Control Systems Magazine}, vol.~42, no.~4, pp. 87--102, 2022.

\bibitem{xu2022efficient}
G.~Xu, G.~Chen, H.~Qi, and Y.~Hong, ``Efficient algorithm for approximating
  {N}ash equilibrium of distributed aggregative games,'' \emph{IEEE
  {T}ransactions on {C}ybernetics}, vol.~53, no.~7, pp. 4375--4387, 2023.

\bibitem{ke2017distributed}
M.~Ke, Y.~Xu, A.~Anpalagan, D.~Liu, and Y.~Zhang, ``Distributed {TOA}-based
  positioning in wireless sensor networks: A potential game approach,''
  \emph{IEEE Communications Letters}, vol.~22, no.~2, pp. 316--319, 2017.

\bibitem{wang2006further}
Z.~Wang, S.~Zheng, S.~Boyd, and Y.~Ye, ``Further relaxations of the {SDP}
  approach to sensor network localization,'' \emph{Tech. Rep.}, 2006.

\bibitem{tseng2007second}
P.~Tseng, ``Second-order cone programming relaxation of sensor network
  localization,'' \emph{SIAM Journal on Optimization}, vol.~18, no.~1, pp.
  156--185, 2007.

\bibitem{marden2009cooperative}
J.~R. Marden, G.~Arslan, and J.~S. Shamma, ``Cooperative control and potential
  games,'' \emph{IEEE Transactions on Systems, Man, and Cybernetics, Part B
  (Cybernetics)}, vol.~39, no.~6, pp. 1393--1407, 2009.

\bibitem{monderer1996potential}
D.~Monderer and L.~S. Shapley, ``Potential games,'' \emph{Games and {E}conomic
  {B}ehavior}, vol.~14, no.~1, pp. 124--143, 1996.

\bibitem{xu2024global}
G.~Xu, G.~Chen, Y.~Hong, B.~Fidan, T.~Parisini, and K.~H. Johansson, ``Global
  solution to sensor network localization: A non-convex potential game approach
  and its distributed implementation,'' \emph{arXiv preprint arXiv:2401.02471},
  2024.

\bibitem{jin2020local}
C.~Jin, P.~Netrapalli, and M.~Jordan, ``What is local optimality in
  nonconvex-nonconcave minimax optimization?'' in \emph{International
  Conference on Machine Learning}.\hskip 1em plus 0.5em minus 0.4em\relax PMLR,
  2020, pp. 4880--4889.

\bibitem{heusel2017gans}
M.~Heusel, H.~Ramsauer, T.~Unterthiner, B.~Nessler, and S.~Hochreiter, ``Gans
  trained by a two time-scale update rule converge to a local {N}ash
  equilibrium,'' \emph{{A}dvances in {N}eural {I}nformation {P}rocessing
  {S}ystems}, vol.~30, 2017.

\bibitem{pang2011nonconvex}
J.-S. Pang and G.~Scutari, ``Nonconvex games with side constraints,''
  \emph{SIAM Journal on Optimization}, vol.~21, no.~4, pp. 1491--1522, 2011.

\bibitem{anderson2010formal}
B.~D. Anderson, I.~Shames, G.~Mao, and B.~Fidan, ``Formal theory of noisy
  sensor network localization,'' \emph{SIAM Journal on Discrete Mathematics},
  vol.~24, no.~2, pp. 684--698, 2010.

\bibitem{fidan2010closing}
B.~Fidan, J.~M. Hendrickx, and B.~D. Anderson, ``Closing ranks in rigid
  multi-agent formations using edge contraction,'' \emph{International Journal
  of Robust and Nonlinear Control}, vol.~20, no.~18, pp. 2077--2092, 2010.

\bibitem{tay1985generating}
T.-S. Tay and W.~Whiteley, ``Generating isostatic frameworks,''
  \emph{Structural Topology 1985 N{\'u}m 11}, 1985.

\bibitem{chen2021distributed}
G.~Chen, Y.~Ming, Y.~Hong, and P.~Yi, ``Distributed algorithm for
  $\varepsilon$-generalized {N}ash equilibria with uncertain coupled
  constraints,'' \emph{Automatica}, vol. 123, p. 109313, 2021.

\bibitem{calafiore2010distributed}
G.~C. Calafiore, L.~Carlone, and M.~Wei, ``Distributed optimization techniques
  for range localization in networked systems,'' in \emph{49th IEEE Conference
  on Decision and Control (CDC)}.\hskip 1em plus 0.5em minus 0.4em\relax IEEE,
  2010, pp. 2221--2226.

\bibitem{shi2010distributed}
Q.~Shi, C.~He, H.~Chen, and L.~Jiang, ``Distributed wireless sensor network
  localization via sequential greedy optimization algorithm,'' \emph{IEEE
  Transactions on Signal Processing}, vol.~58, no.~6, pp. 3328--3340, 2010.

\bibitem{anderson2008rigid}
B.~D. Anderson, C.~Yu, B.~Fidan, and J.~M. Hendrickx, ``Rigid graph control
  architectures for autonomous formations,'' \emph{IEEE Control Systems
  Magazine}, vol.~28, no.~6, pp. 48--63, 2008.

\bibitem{cao2021bearing}
K.~Cao, Z.~Han, Z.~Lin, and L.~Xie, ``Bearing-only distributed localization: A
  unified barycentric approach,'' \emph{Automatica}, vol. 133, p. 109834, 2021.

\bibitem{facchinei2010penalty}
F.~Facchinei and C.~Kanzow, ``Penalty methods for the solution of generalized
  {N}ash equilibrium problems,'' \emph{SIAM Journal on Optimization}, vol.~20,
  no.~5, pp. 2228--2253, 2010.

\bibitem{koshal2016distributed}
J.~Koshal, A.~Nedi{\'c}, and U.~V. Shanbhag, ``Distributed algorithms for
  aggregative games on graphs,'' \emph{Operations Research}, vol.~64, no.~3,
  pp. 680--704, 2016.

\bibitem{gao2017canonical}
D.~Y. Gao, V.~Latorre, and N.~Ruan, \emph{Canonical Duality Theory: Unified
  Methodology for Multidisciplinary Study}.\hskip 1em plus 0.5em minus
  0.4em\relax Springer, 2017.

\bibitem{xu2023algorithm}
G.~Xu, G.~Chen, and H.~Qi, ``Algorithm design and approximation analysis on
  distributed robust game,'' \emph{Journal of Systems Science and Complexity},
  vol.~36, no.~2, pp. 480--499, 2023.

\bibitem{10044227}
G.~Chen, G.~Xu, W.~Li, and Y.~Hong, ``Distributed mirror descent algorithm with
  {B}regman damping for nonsmooth constrained optimization,'' \emph{IEEE
  Transactions on Automatic Control}, vol.~68, no.~11, pp. 6921--6928, 2023.

\bibitem{23}
D.~Dua and C.~Graff,
  \url{http://archive.ics.uci.edu/dataset/343/ujiindoorloc+m\\ ag}, Accessed on
  2019.

\end{thebibliography}

\end{document}